\documentclass[twocolumn,a4paper]{jjspeme}
\usepackage{graphicx}
\usepackage{amsmath,amssymb,amsthm}
\usepackage{times}
\usepackage{type1cm}
\usepackage{float,fancyhdr}

\theoremstyle{definition}

\newcommand{\R}{\mathbb{R}}


\begin{document}

\setlength{\baselineskip}{13.0pt}

\title{A closed linkage mechanism having the shape of a discrete M\"obius strip}
\etitle{}
  
\author{
Shizuo KAJI${}^{*}$
}

\eauthor{
}


\begin{abstract}
A closed linkage mechanism in three-dimensional space is an object comprising rigid bodies connected with hinges in a circular form like a rosary. Such linkages include Bricard6R and Bennett4R. To design such a closed linkage, it is necessary to solve a high-degree algebraic equation, which is generally difficult. In this lecture, the author proposes a new family of closed linkage mechanisms with an arbitrary number of hinges as an extension of a certain Bricard6R. They have singular properties, such as one-dimensional degree of freedom (1-DOF), and certain energies taking a constant value regardless of the state. These linkage mechanisms can be regarded as discrete M\"obius strips and may be of interest in the context of pure mathematics as well. However, many of the properties described here have been confirmed only numerically, with no rigorous mathematical proof, and should be interpreted with caution.
\end{abstract}

\maketitle

\begingroup  \makeatletter
    \def\@makefntext#1{\parindent=\z@#1}%
    \footnotetext{%
\setlength{\tabcolsep}{1pt}
\begin{tabular}{rrll}
${}^*$ &Yamaguchi University& / JST, PRESTO \\
\end{tabular}\\
    }%
\endgroup

\section{Introduction}
\begin{figure}[!htb]
\begin{center}
\includegraphics[width=0.9truein,clip]{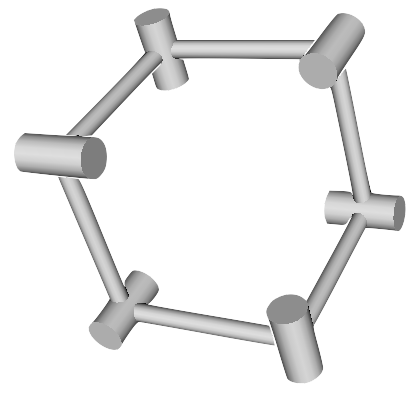}\hspace{1cm}
\includegraphics[width=0.9truein,clip]{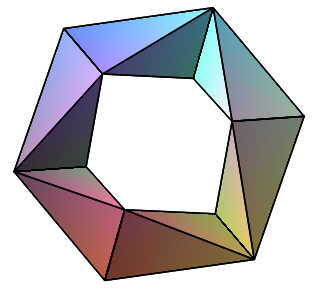}
\end{center}
\caption{Bricard6R and $6$-Kaleidocycle}
\label{fig:opensim-K6}
\end{figure}

\headheight 10pt              
\headsep 5pt                 
\thispagestyle{fancy}
  \rhead{
  \small Translated from the Symposium Proceedings of the 2018 Spring meeting of the Japan Society for Precision Engineering,
  pp. 62--65, 1 Mar. 2018.}
A joint is called a \emph{hinge} (or revolute joint) if it has a rotation axis and can freely change its direction around this axis. Let us call a set of rigid bodies connected by hinges a \emph{linkage mechanism}.
 One example is the Bricard6R linkage (Fig. \ref{fig:opensim-K6}, left), a linkage consisting of six hinges. One of the most important subjects of the study of linkage mechanisms is to determine all possible states, i.e., the configuration space of the linkage mechanism. Mathematically, because the state of a hinge can be identified with the unit circle\footnote{More precisely, $S^1$ is obtained by identifying the two ends $\{0, 2\pi\}$
 of the closed space $[0, 2\pi]$.} 
$S^1 := [0,2\pi)$
 using the rotation angle, the configuration space of a system consisting of $n$ hinges will be a subspace of the torus $(S^1)^n$.

Linkage mechanisms have various applications, such as in robot arms, folding and deployable mechanisms, and power transmission. Details of the classification and analysis of linkage mechanisms can be found in \cite{Chen2011}.
 In the context of mathematics, studying the configuration space of linkage mechanisms is a topic in the area of topology \cite{Milgram2004,Pavesic2017}, which is the specialty of the author.

The most important property of a configuration space, which should be considered first, is its dimension. This is equal to the degree of freedom (\textbf{DOF}) of the motion of the mechanism\footnote{A configuration space can have a singularity. When we discuss the dimension of a configuration space, however, we will neglect them. We will consider the degrees of freedom assuming the mechanism is in generic states that cover almost all possible states.}.
The \emph{mobility formula} (also called the Chebychev-Gr\"ubler-Kutzbach criterion)
\begin{equation}\label{eq:mobility}
M=6(N-1-n)+\sum_{i=1}^n f_i
\end{equation}
 gives $M$, an
 estimation of the DOF of the system, by simply counting the number of variables and constraints involved. Here, $N$ is the number of rigid bodies involved, $n$ is the number of joints, and $f_i$ is the DOF of each joint\footnote{For a hinge, $f_i=1$}.
 If some of the constraints are redundant and not independent, the system possesses degrees of freedom larger than this estimation; such a system is called \emph{over-constrained}. Some famous mechanisms, including the Bricard6R and Bennett4R linkages, are over-constrained. In contrast, other systems have degrees of freedom smaller than this simple estimation. This occurs, as discussed later, as a result of the system being a real algebraic variety. In this lecture, we will consider closed linkage mechanisms, each of which has $n$ hinges for an arbitrary natural number $n \ge 6$ that are connected in a cyclic manner, and has just one DOF. In particular, the $n = 6$ case corresponds to a variation of the Bricard6R linkage, and hence, the $n\ge 7$ case can be regarded as its generalisation.

\section{Kaleidocycle}
The cyclic, closed linkage mechanism described here is based on a toy called a Kaleidocycle \cite{Schattschneider1987}({\figurename\ref{fig:opensim-K6}}, right). 
A Kaleidocycle\footnote{A Kaleidocycle is also referred to as a rotating ring of tetrahedra in some articles; be careful when performing a literature search.}, which can be made by folding a sheet of paper, is a cyclic object consisting of $n$ congruent tetrahedrons, where two opposite sides of each tetrahedron serve as hinges to connect the adjacent tetrahedrons, and you can continue rotating the object\footnote{Hereinafter, this motion will be referred to as rotating motion; it might also be called everting motion.}. 
Let us call this linkage with $n$ hinges for $n \ge 6$ the $n$-Kaleidocycle.
The main subject of this lecture is to show that the DOF of a Kaleidocycle can be made one-dimensional (rotating motion) by selecting a particular shape for the tetrahedron, 
although a Kaleidocycle generally has a larger DOF, such as bending motion in addition to rotating, as $n$ becomes larger. We first describe a mathematical framework for Kaleidocycles using the Denavit-Hartenberg parameters, which are the standard parameters for the analysis of linkage mechanisms (e.g., \cite[Ex 5.2, Ex 8.13]{Bertini}).

We first consider a serial linkage mechanism obtained when opening a closed linkage by cutting one of its hinges, and we will treat conditions for the serial linkage to be closed as an inverse kinematics problem. Let the direction of the $i$-th $(0\le i\le n)$ hinge be $b_i \in S^2 := \{(x, y, z) \in \R^3 | x^2 + y^2 + z^2 = 1\}$ and the position of its centre be $\gamma_i\in \R^3$. Here, the $n$-th hinge corresponds to the end that emerges as a result of opening up a ring. The serial linkage considered will close if the opposite end overlaps with this, and there are two kinds of those closed cases. Let us say that such a closed linkage is
\begin{itemize}  \setlength{\itemsep}{-0.5em}
\item \textbf{Oriented} if $b_{n}=b_0$, $\gamma_{n}=\gamma_{0}$
\item \textbf{Non-oriented} if $b_{n}=-b_0$, $\gamma_{n}=\gamma_{0}$.
\end{itemize}
Because, for the linkage mechanism we are now considering ($n$-Kaleidocycle), two adjacent hinges are fixed with a rigid body, and the twist angle for two adjacent hinges is constant across all sets of two adjacent hinges (as all the tetrahedrons are congruent), we can write $b_{i-1}\cdot b_i=c$, using a parameter $c \in [-1, 1]$ that stands for the cosine of the twist angle, and the inner product of $b_{i-1}$ and $b_i$. 
If $b_{i-1}$ and $b_i$ are not parallel, $\gamma_i$ can inductively be determined as $\gamma_{i-1}+b_{i-1}\times b_i$ using the cross product\footnote{For the purpose of discussing the DOF of a system alone, the center of each hinge does not have to be physically positioned at $\gamma_i$, and the hinge is allowed to be at any position in the line that passes through $\gamma_i$ and possesses orientation vector $b_i$ (see Fig. \ref{fig:extremeK}, right). This is an important remark when considering applications to folding structures for which each hinge can slide along its axis.}. If all hinges are parallel (i.e., $c = \pm 1$), the corresponding linkage practically represents a planar one, which may not be of interest; such a case will be called a trivial case and neglected hereafter.

Summarising the above, the configuration space of an $n$-Kaleidocycle can be identified with the solution space (set of real solutions) of the simultaneous quadratic equations
 \begin{equation}\label{eq:kaleidocycle}
\sum_{i=1}^n b_{i-1}\times b_i=(0,0,0), \quad b_{i-1}\cdot b_i=c, \quad  b_{i}\cdot b_i=1 \ (1\le i \le n).
\end{equation}
Contrary to their appearance, studying this solution space is not easy (cf. \cite{Raghavan1993}).

It is worth noting that the parameter $c$ cannot be freely chosen. For some tetrahedrons, it is impossible to connect copies of them in a way such that both ends meet to make a cyclic object. For example, when $n = 6$, the above simultaneous equations have solutions only when 
$c = \pm 1$ for oriented cases, or $c = 0$ for non-oriented cases. The non-oriented case with $c = 0$ corresponds to the Bricard6R linkage and has a configuration space homeomorphic to $S^1$, i.e., this linkage has one degree of freedom. For each $n$ and each of the oriented and non-oriented cases, the range of $c$ within which solutions exist can be written as follows, using a certain constant $c_n$:
\begin{itemize}     \setlength{\itemsep}{-0.5em}
\item Oriented cases with even $n$:  $[-1,1]$
\item Oriented cases with odd $n$: $[-c_n,1]$
\item Non-oriented cases with even $n$: $[-c_n,c_n]$
\item Non-oriented cases with odd $n$: $[-1,c_n]$
\end{itemize}
$n$-Kaleidocycles with a non-trivial boundary value $c = \pm c_n$ in the above range are called \emph{extreme} $n$-Kaleidocycles, which are the main topic of this lecture. The oriented and non-oriented extreme cases with an odd value of $n$ are transformed to each other if every other hinge in one linkage is oriented toward the opposite direction; the $c = \pm c_n$ cases with an even value of $n$ are also transformed to each other if every other hinge in one linkage is oriented toward the opposite direction. Therefore, there is essentially only one\footnote{There also exists the mirror for a $c$.} extreme $n$-Kaleidocycle for each $n$.

\begin{figure}[!htb]
\begin{center}
\includegraphics[width=0.8truein,clip]{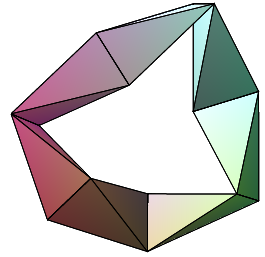}\hspace{1cm}
\includegraphics[width=0.8truein,clip]{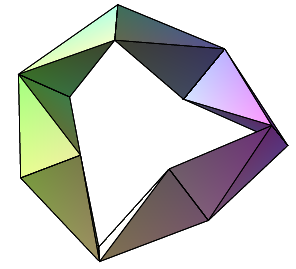}
\end{center}
\caption{Oriented (left) and non-oriented (right) extreme 7-Kaleidocycles}
\end{figure}

The value of $c_n$, which is a solution of a high-degree polynomial equation, cannot be obtained algebraically\footnote{When $n = 7$, $c_n$ can be determined algebraically, as it is a solution of a cubic equation.}, 
but it can be easily obtained numerically as a solution of a constrained optimisation problem. Values of $c_n$ obtained by numerical computation are shown in {\bf \tablename\ref{table:tabtest}}.

\section{extreme Kaleidocycle}
In the previous section, we saw that the set of all possible states of $n$-Kaleidocycles can be identified with the solution space of the simultaneous quadratic equations \eqref{eq:kaleidocycle}. Considering the rotational and reflectional symmetries of the entire system, we can choose $b_0=(0,0,1)$ and $b_1=(0,\sqrt{1-c^2},c)$. With this, let us count the DOF. Because the number of variables (except $c$) is $3(n - 2)$, and the number of equations is $3+(n-1)+(n-2)$, the DOF of the system is expected to be $n-6$, which is consistent with \eqref{eq:mobility}. Although the expected DOF for a $6$-Kaleidocycle with $c = 0$ is $0$, the actual DOF is $1$, and this system is an over-constrained one (to account for this, \cite{Foweler-Guest2005} explores improvements for Eq. \eqref{eq:mobility}). The phenomenon where the expected DOF differs from the actual DOF occurs as a singular property for the extreme $n$-Kaleidocycle, and the DOF is in fact $1$, regardless of $n$.
 $n$-Kaleidocycles other than the extreme ones can undergo bending motions in addition to the rotating motion, and generically have 
 $(n - 6)$-DOF, according to \eqref{eq:mobility}.
  Although the author is not very knowledgeable about linkage mechanisms, he presumes that the extreme $n$-Kaleidocycles might represent a very rare 
family of linkage mechanisms, as he could not find another example of a mechanism for which the DOF is always $1$ no matter how many hinges the linkage has.

\textbf{Remark:  }
The set of all possible Kaleidocycles can be identified with the space of all real solutions for the simultaneous equations \eqref{eq:kaleidocycle} with the parameter $c$ being a variable. For each fixed value of $c$, a slice is determined in this solution space, and the connected component of this slice provides the configuration space of the corresponding Kaleidocycle. The configuration space of the extreme $n$-Kaleidocycle corresponds to the slice where $c$ has the boundary value, and the slice determined by $c$ reduces to a one-dimensional space only there\footnote{Thinking of the real algebraic variety determined by \eqref{eq:kaleidocycle} as the moduli space for a Kaleidocycle, you are seeing a singular fibre for the projection map onto the $c$ axis.}. 
Because this is an important point, let us illustrate the situation with a very simple example by devoting some pages.
For simplicity, in this example only, we consider joints, not hinges, that can orient in any direction (i.e., have the same DOF as $S^2$).

\begin{figure}[!htb]
\begin{center}
\includegraphics[width=1.3truein,clip]{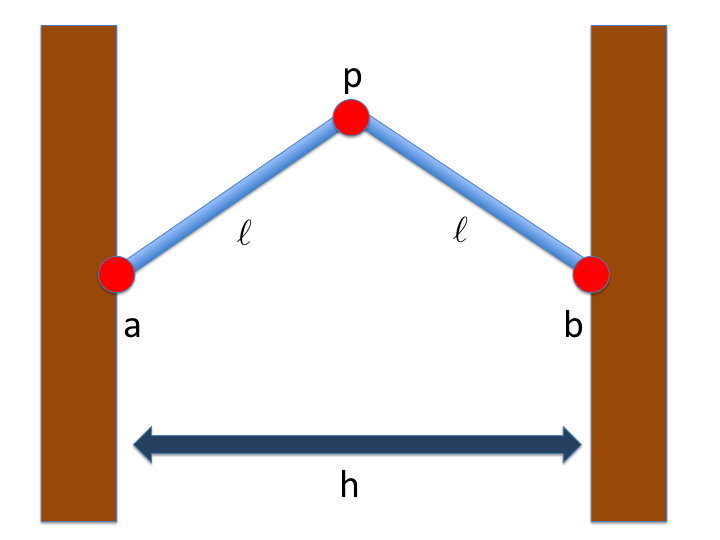}
\end{center}
\caption{A simple example in which the DOF change according to parameter values.}
\label{fig:wall}
\end{figure}

As shown in Fig. \ref{fig:wall}, we consider a system consisting of three joints and two rods, in which each of the two rods is connected to a wall with a joint and the opposite ends of the two rods are connected to each other with another joint $p$
 in the three-dimensional space\footnote{Equation \eqref{eq:mobility} cannot be applied to this system without modifications because the walls are fixed.}. The configuration space of this system can be identified with the space of all possible locations of the junction $p$, which is equal to the space of solutions of the simultaneous quadratic equations $\{p\in \R^3\mid |p-a|^2=l^2, |p-b|^2=l^2\}$.
  Regarding this space, there are three possible cases, according to the length of the rod $l$:
\begin{itemize}     \setlength{\itemsep}{-0.5em}
\item When $2l > h$, the configuration space is a circle (in the plane normal to the page).
\item When $2l = h$, the configuration space is a point.
\item When $2l < h$, the configuration space is an empty set.
\end{itemize}
Note that the dimension of the configuration space reduces when the parameter $l$ has the boundary value $h/2$. Real solutions exhibit more complicated behaviours compared with complex solutions of a system of algebraic equations.

Let us describe some of the intriguing properties of an extreme $n$-Kaleidocycle, other than the DOF being one. In \cite{Safsten2016}, a system was considered in which a torsional spring that has an energy potential proportional to the square of its rotation angle is attached to each hinge of a $6$-Kaleidocycle. One subject studied in the reference was to determine
equilibrium points for this $6$-Kaleidocycle because the energy of this linkage changes according to its state during its rotating motion. Meanwhile, an extreme $n$-Kaleidocycle with seven or more hinges, in contrast, has a peculiar property: the energy is constant, i.e., when $n \ge 7$,
\[
 E_{bend} :=  \sum_{i=1}^n \arccos\left(\dfrac{(b_{i-1}\times b_i)\cdot (b_i\times b_{i+1})}
 {|b_{i-1}\times b_i||b_i\times b_{i+1}|} \right)^2
\]
is almost constant in the configuration space of an extreme $n$-Kaleidocycle. This is interesting with respect to both theory and application, as it means that although the angle of each hinge will have various values as the Kaleidocycle rotates, the total energy potential remains constant, and no force is required to rotate the linkage.

An extreme n-Kaleidocycle ($n \ge 7$) has a constant value also for the following two energies in its configuration space. One is the Coulomb potential, when the centre of each hinge is electrically charged:
\[
E_{clmb} :=\sum_{i<j} \dfrac{1}{|\gamma_i-\gamma_j|^\alpha} \qquad \alpha\in \R
\]
The other is the energy potential when an oriented extreme $n$-Kaleidocycle has a very small dipole at the centre of each hinge along the direction of the hinge:
\[
E_{dipl} :=\sum_{i<j} \dfrac{b_i\cdot b_j}{|\gamma_i-\gamma_j|^3} - \dfrac{3(b_i\cdot (\gamma_i-\gamma_j))(b_j\cdot (\gamma_i-\gamma_j))}{|\gamma_i-\gamma_j|^5}
\]
It seems more mysterious that the total energy is constant for these two systems, considering the fact that interactions between any two hinges contribute to the energy in this system, in contrast with the above-mentioned system with torsional springs, in which the energy depends only on the interactions between neighbouring hinges.

Table \ref{table:tabtest} shows the values of $c_n$, twist (described later; when $c = c_n$), 
$E_{bend}$, and $E_{dipl}$ (oriented cases) for some $n \ge 6$.
\begin{table}[htb!]
\fontsize{7.5pt}{11pt}\selectfont
\caption{Some key values for extreme $n$-Kaleidocycles}
\begin{center}
\begin{tabular}{|c|c|c|c|c|c|c|} \hline
$n$ & 6 & 7 & 8 & 9 & 15 & 38 \\ \hline
$c_n$ & 0 & 0.2954 & 0.4700 & 0.5852 & 0.8533 & 0.9773 \\ \hline
$Tw$ & 1.500 & 1.416 & 1.377 & 1.355 & 1.309 & 1.291 \\ \hline
$E_{bend}$ & varies & 11.9 & 10.4 & 9.24 & 5.60 & 2.23\\ \hline
$E_{dipl}$ & NA & -4.23 & NA & -10.0 & -83.7 & NA \\ \hline
\end{tabular}
\end{center}
\label{table:tabtest}
\end{table}

If considering all Kaleidocycles, including non-extreme cases, $E_{bend}$ appears to have a minimum value when the positions of $\gamma_i$ are as planar as possible and have a higher symmetry\footnote{For example, when $n$ is even, $E_{bend}$ has a minimum value for the Kaleidocycle with $c = 0$ in which every other hinge is parallel to the $z$-axis and they form a regular polygon when viewed from the direction of the $z$-axis.}, and an extreme $n$-Kaleidocycle (except when $n = 7$) does not provide even a local minimum.

\section{Discrete Strips}
As inferred in the title of this article, extreme $n$-Kaleidocycles have some geometrical aspects, which are the topic of this last section. The author apologises for the description being a little bit rough, partly because of the limited number of pages.
Consider a closed curve $\gamma:S^1\to \R^3, |\dot\gamma(s)|=1$ with a length of $2\pi$ in the three-dimensional space. 
When a map $b: [0,2\pi] \to S^2$ such that $\forall s, \dot{\gamma}(s)\cdot b(s)=0$ is given, we say that the set of $\gamma$ and $b$ is a \emph{strip with $\gamma$ being its centre line}. Specifically, a strip is a curve for which a direction normal to the tangent has been determined at each point on the curve. Imagine a M\"obius strip-like object that is made by gluing both ends of a long and narrow sheet of paper together after twisting it several times.

For a given strip, we consider the \emph{number of half-twists} for the strip; although its meaning may be intuitively clear, it is precisely determined as $2(Tw + Wr)$ using the C\v alug\v areanu-White theorem. Here,
\begin{align*}
Tw & :=\frac{1}{2\pi}\int_0^{2\pi} \dot{b}(s)\cdot (\dot{\gamma}(s)\times b(s)) ds \\
Wr & :=\frac{1}{4\pi}\int_0^{2\pi} \int_0^{2\pi} \frac{(\dot{\gamma}(s_1)\times \dot{\gamma}(s_2))\cdot (\dot{\gamma}(s_1)-\dot{\gamma}(s_2))}
{|\dot{\gamma}(s_1)-\dot{\gamma}(s_2)|^3}ds_1 ds_2
\end{align*}
and these quantities are called the twist and writhe, respectively. Each of these quantities has various geometrical interpretations. For example, $Tw$ measures how much the vector obtained by parallel-transporting $b(0)$ along the centre line one round differs from $b(2\pi)$, while $Wr$ represents the average number of self-intersections in the image of the centre line projected in every direction. In particular, the latter is a quantity that does not depend on $b$, only on the centre line. Intuitively, $Tw$ and $Wr$ represents the extent of the strip winding around the centre line and the extent of the centre line running irregularly, respectively, and the sum makes the twisting number for the whole system.

\begin{figure}[!htb]
\begin{center}
\includegraphics[width=1.0truein,clip]{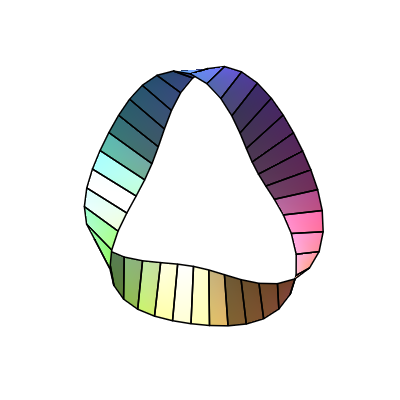}
\includegraphics[width=1.0truein,clip]{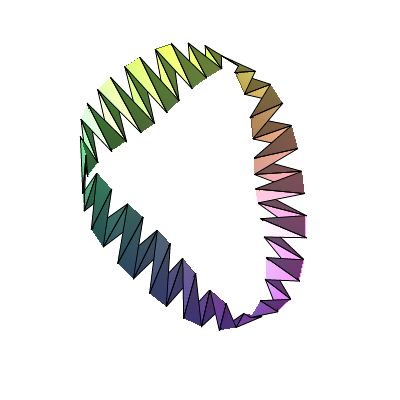}
\includegraphics[width=1.0truein,clip]{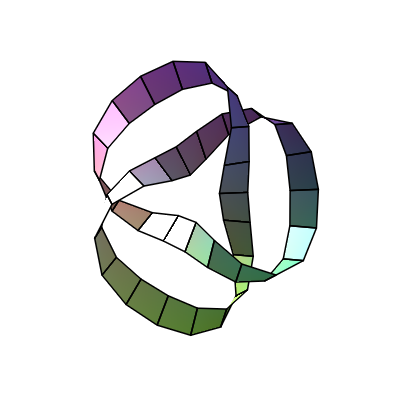}
\end{center}
\caption{
Strip representations for the extreme $38$-Kaleidocycles corresponding to the $c = c_n$ and $c = -c_n$ cases (left and centre, respectively). The $38$-Kaleidocycle (right) made by sliding the centre of a hinge of the extreme $19$-Kaleidocycle along the hinge with a certain offset and extending it one more round can be regarded as a M\"obius strip cut along its centre line, which gives a rotating trifoliate knot (cf. \cite{trefoil}).}
\label{fig:extremeK}
\end{figure}

A Kaleidocycle can be regarded as a strip such that its centre line is the polygonal line connecting $\gamma_i$ and there are associated vectors interpolating $b_i$ and $b_{i + 1}$ with a constant angular velocity on the line segment $\gamma_i\gamma_{i+1}$ (Fig. \ref{fig:extremeK}, left). For a continuous curve, the direction in which the tangent changes is called the normal, and the cross product of the tangent and normal is called the \emph{binormal}; as $b_i$ is normal to both of the two neighbouring ``tangents'' $\gamma_{i-1}\gamma_i, \gamma_i\gamma_{i+1}$, it can be thought of as the discrete binormal up to sign. Furthermore, the rate of change of a binormal is called the \emph{torsion}; $\pm c$ can be thought of as its discrete version up to a constant factor.

In this setting, calculation shows 
$Tw = \frac{1}{2\pi}\sum_{i=1}^n \arccos(b_{i-1} \cdot b_{i})=\frac{n\arccos(c)}{2\pi}$,
 and $Wr$ for the central polygonal line is given, e.g., in \cite{Klenin}. For an extreme $n$-Kaleidocycle, the number of half-twists is $3$ for $c = c_n$ or $n - 3$ for $c = -c_n$. The $c = c_n$ case gives, by definition, the least twisted strip among those with an odd number of half-twists, and it seems mysterious that the number of half-twists for this case is not $1$, but $3$. As seen from Table \ref{table:tabtest}, $Tw$ for the $c = c_n$ case appears to monotonically converge to a constant as $n$ increases; if this is the case, $2(Tw + Wr) = 3$ requires that $Wr$ also monotonically converges to a constant.
It is unknown what the strip obtained in the limit and its centre line are. With its rotating motion, its centre line also moves, but because $Tw$ is constant, $Wr$ is also constant. It is interesting to see this from the point of view of the curve on the sphere drawn by the tangent moving on the centre line, i.e., the Gauss map $s\mapsto \dot{\gamma(s)}$. Here, the discrete Gauss map is supposed to be yielded by the curve connecting the vectors on the sphere with arcs of great circles. From the Gauss--Bonnet theorem, the area cut off by this curve from the sphere on the left hand is equal to $2\pi Wr$ modulo $2\pi$. $E_{bend}$ is equal to the squared sum of the length of each arc. According to the rotating motion, the image of the Gauss map moves in such a way that both the area cut off and the squared sum of the length of each arc are constant. Although the author hasn't been able to write down the time evolution of this yet\footnote{The rotating motion takes the form of so-called falling cat motion, the motion resulting in a different direction as a whole after rotating once while preserving the angular momentum. Therefore, applications to, e.g., antennas that can change their direction under the condition of zero gravity might be possible.},
 results of numerical experiments indicate soliton-like structures, which are expected to provide a key to understanding the constancy of the various energies.

The motion of the centre line can be interpreted also in the context of an elastic rod problem. If we consider the tangent and $b$ and their cross product at each point on the centre line of a strip, we have a family of orthonormal systems (\textbf{adapted frame}) along the centre line. A curve with this adapted frame is often studied as a model of an elastic rod (cf. \cite{Langer1996}). For an extreme $n$-Kaleidocycle with $n \ge 7$, as both $b_{i-1} \cdot b_{i}$ and $E_{bend}$ are constant, we can see that the isotropic Kirchhoff elastic energy also remains constant during the rotating motion. Implications for the theory in \cite{Langer1996} may also be of much interest.

Finally, because the number of half-twists is always an integer, two solutions of \eqref{eq:kaleidocycle} with different numbers of half-twists belong to different connected components. The properties of those connected components may be of mathematical interest; the author numerically searched for $n$-Kaleidocycles with $c$ having a boundary value on the components other than those with the number of half-twists being $3$ or $n - 3$ and found that at least they seem to have no singular properties, similar to those for the extreme cases described above.

\begin{figure*}[!h]
 \begin{center}
  \includegraphics[height=18.3cm]{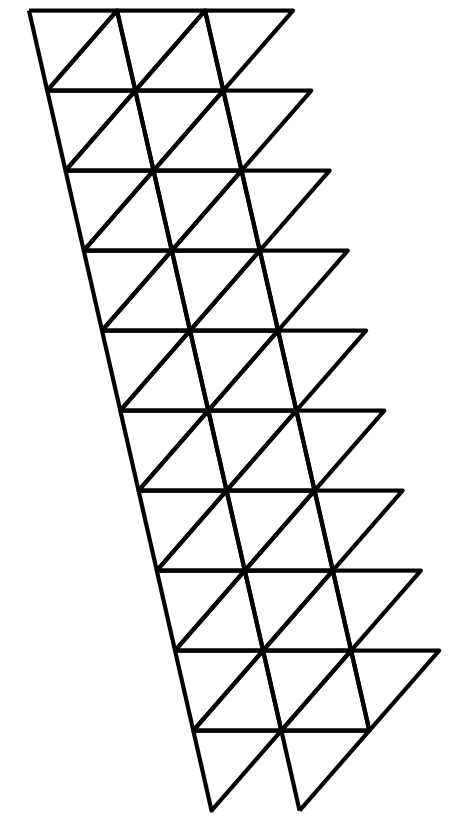}
  \caption{Development of extreme $9$-Kaleidocycle. The right and lower ends are margins for gluing.}
 \end{center}
\end{figure*}

\section*{Acknowledgements}
This lecture is partially based on a collaborative research with E. Fried, M. Grunwald, and J. Sch\"onke (OIST). 
The author is supported by JST PRESTO (Grant Number: JPMJPR16E3).



\begin{thebibliography}{9}

\bibitem{Foweler-Guest2005}
P. W. Fowler and S. D. Guest,
{\it A symmetry analysis of mechanisms in rotating rings of tetrahedra},
Proc. R. Soc. A 461, 1829--1846 (2005).

\bibitem{Klenin}
K. Klenin and J. Langowski,
{\it Computation of writhe in modeling of supercoiled DNA},
 Biopolymers 54, 307--317 (2000).

\bibitem{Langer1996}
J. Langer and D. Singer,
{\it Lagrangian Aspects of the Kirchhoff Elastic Rod},
SIAM Reviews 38, 605--618 (1996).

\bibitem{Milgram2004}
R. J. Milgram and J. C. Trinkle, 
{\it The geometry of configuration spaces for closed chains in two and three dimensions},
Homology, Homotopy and Applications 6.1, 237--267 (2004).

\bibitem{Pavesic2017}
P. Pave\v si\'c,
{\it A Topologist's View of Kinematic Maps and Manipulation Complexity},
arXiv:1707.03899.

\bibitem{Raghavan1993}
M. Raghavan and B. Roth,
{\it Inverse Kinematics of the General 6R Manipulator and Related Linkages},
J. Mech. Des 115(3), 502--508 (1993).


\bibitem{Safsten2016}
C. Safsten, T. Fillmore, A. Logan, D. Halverson and L. Howell
{\it Analyzing the Stability Properties of Kaleidocycles},
J. Appl. Mech 83(5), 051001 (2016).

\bibitem{Schattschneider1987}
D. Schattschneider and W. M. Walker, 
{\it M. C. Escher Kaleidocycles}, Pomegranate Communications: Rohnert Park, CA (1987).
(TASCHEN; Reprint edition (2015)).
%

\bibitem{Bertini}
A. J. Sommese, J. D. Hauenstein, D. J. Bates, and C. W. Wampler,
{\it Numerically Solving Polynomial Systems with Bertini},
Software, Environments, and Tools, Vol. 25, SIAM, Philadelphia, PA (2013).
%

\bibitem{Chen2011}
Z. You and Y. Chen,
{\it Motion Structures: Deployable Structural Assemblies of Mechanisms},
Taylor \& Francis 2011.

\bibitem{trefoil}
D. M. Walba et al.,
{\it The thyme polyethers : An approach to the synthesis of a molecular knotted ring},
Tetrahedron 42, 1883--1894 (1986).

\end{thebibliography}
\end{document}